\documentclass[11pt]{article}

\usepackage{latexsym,mathrsfs}
\usepackage{amsmath,amssymb}
\usepackage{amsthm,enumerate,verbatim}
\usepackage{amsfonts}
\usepackage{graphicx}
\usepackage{algorithm}
\usepackage{algorithmic}
\usepackage{url}

\newtheorem{theorem}{Theorem}

\newtheorem{definition}{Definition}
\newtheorem*{definition*}{Definition}

\DeclareMathOperator{\conv}{conv} 
 
\DeclareMathOperator{\argmin}{argmin} 
 
\DeclareMathOperator{\rank}{rank}

\DeclareMathOperator{\col}{col}

\DeclareMathOperator{\tr}{tr} 
\DeclareMathOperator{\diag}{diag}

\title{Introduction to Nonnegative Matrix Factorization} 
\date{}

\author{Nicolas Gillis \\ 
Department of Mathematics and Operational Research \\ 
Facult\'e Polytechnique, Universit\'e de Mons \\ 
Rue de Houdain 9, 7000 Mons, Belgium\\
 nicolas.gillis@umons.ac.be  
}

\begin{document}

\maketitle

\begin{abstract} 
In this paper, we introduce and provide a short overview of nonnegative matrix factorization (NMF). Several aspects of NMF are discussed, namely, the application in hyperspectral imaging, geometry and uniqueness of NMF solutions, complexity, algorithms, and its link with extended formulations of polyhedra. In order to put NMF into perspective, the more general problem class of constrained low-rank matrix approximation problems is first briefly introduced.   
\end{abstract}

\section{Introduction} 

Constrained low-rank matrix approximation (CLRMA) is becoming more and more 
popular because it is able to extract pertinent information from 
large data sets; 
see, for example, the recent survey~\cite{UHZB14}. 
CLRMA is equivalent to \emph{linear dimensionality reduction}. 
Given a set of $n$ data points $m_j \in \mathbb{R}^{p}$ ($j=1,2,\dots,n$), the 
goal is to find a set of $r$ basis
vectors $u_k \in \mathbb{R}^{p}$ ($k=1,2,\dots,r$) and the corresponding weights 
$v_{kj}$ so that for all $j$, 
$
m_j \; \approx  \; \sum_{k = 1}^r v_{kj} \, u_k . 
$ 
This problem is equivalent to the low-rank approximation of matrix $M$, with 
\[ 
M = [m_1 \; m_2 \; \dots \; m_n] \; \approx \; [u_1 \; u_2 \; \dots \; u_k] 
[v_1 \; v_2 \; \dots \; v_n] = UV, 
\] 
where each column of $M$ is a data point, each column of $U$ is a basis vector, 
 and each column
of $V$ provides the coordinates of the corresponding column of $M$ in the basis 
$U$. 
In other words, each column of $M$ is approximated by a linear combination of 
the columns of $U$. 

In practice, when dealing with such models, two key choices exist:  
\begin{enumerate} 

\item \emph{Measure of the error} $M-UV$. Using the standard least-squares
error, $\|M-UV\|_F^2 = \sum_{i,j}(M-UV)^2_{i,j}$, leads to principal 
component analysis (PCA) that can be
solved by using the singular value decomposition (SVD). 
Surprisingly, one can show that the optimization problem in variables 
$(U,V)$ 
has no spurious local minima (i.e., all local minima are global), which 
explains why it can be solved efficiently despite the error being nonconvex. 
Note that the resulting problem can be reformulated as a semidefinite program 
(SDP)
by using the Ky Fan $2$-$k$-Norm~\cite[Prop.~2.9]{DV16}. 

If data is missing or if weights are assigned to the entries of $M$, the problem 
can be cast as a weighted low-rank matrix approximation (WLRA) 
problem with error $\sum_{i,j} W_{i,j} (M-UV)^2_{i,j}$ for some nonnegative 
weight matrix $W$, where $W_{i,j} = 0$ when the entry $(i,j)$ is 
missing~\cite{SJ03}.  
Note that if $W$ contains entries only in $\{0,1\}$, then the problem is also 
referred to as PCA with missing data or low-rank matrix completion with noise. 

WLRA is widely used for recommender systems~\cite{KBV09} that predict the 
preferences of users for a given product based
on the product's attributes and user preferences. 

If the sum of the absolute values of the entries of the error $\sum_{i,j} 
|M-UV|_{i,j}$ is used, we obtain
yet another variant more robust to outliers (sometimes referred to as robust 
PCA~\cite{CLM11}). It can be used, for example, for background subtraction in 
video sequences where the noise (the moving objects) is assumed to be sparse 
while the background has low rank.   

\item \emph{Constraints that the factors $U$ and $V$ should satisfy}. These 
constraints depend on
the application at hand and allow for meaningful interpretation of the factors. 
For example,
$k$-means\footnote{$k$-means is the problem of finding a set of centroids $u_k$ 
such that the sum of the distances between each data point and the closest 
centroid is minimized.} is equivalent to requiring the factor $V$ to have a 
single nonzero entry in each column that is equal to one, so that the columns 
of $U$ are cluster centroids. 
Another widely used variant is sparse PCA, which requires that the
factors ($U$ and/or $V$) be sparse~\cite{D07, JNR10, LT11}, thus yielding a 
more compact and easily interpretable decomposition (e.g., if $V$ 
is sparse, each data point is the linear combination of only a few basis 
elements). 
Imposing componentwise nonnegativity on both factors $U$ and $V$ leads to
nonnegative matrix factorization (NMF). 
For example, in document analysis where each column of $M$ corresponds to a 
document (a vector of word counts), these nonnegativity constraints allow one to 
interpret the columns of the factor $U$ as topics, 
and the columns of the factor $V$ indicate in which proportion each document
discusses each topic~\cite{LS99}. In this paper, we focus on this 
particular variant of CLRMA. 
\end{enumerate}  

CLRMA problems are at the heart of many fields of applied mathematics and 
computer science, including,   
statistics and data analysis~\cite{J2002}, 
machine learning and data mining~\cite{elden2007matrix}, 
signal and image processing~\cite{Mag14}, 
graph theory~\cite{chung1997spectral}, numerical linear algebra, 
and systems theory and control~\cite{markovsky2011low}. 
The good news for the optimization community is that these CLRMA models lead to 
a wide variety of 
theoretical and algorithmic challenges for optimizers: 
Can we solve these problems? Under which conditions? What is the most 
appropriate model for a given application? 
Which algorithm should we use in which situation? What type of guarantees can 
we 
provide?  

CLRMA problems can be formulated in the following way: 
\begin{equation} \label{clrma}
\min_{U \in \Omega_U, V \in \Omega_V} \|M-UV\| . 
\end{equation}  
As an introduction, below we discuss several aspects of~\eqref{clrma}.  
\newline

\noindent
{\bf Complexity.}
As soon as the norm $\| \cdot \|$ is not the Frobenius norm or the feasible 
domain 
has constraints (i.e., $\Omega_U \neq \mathbb{R}^{p \times r}$ or $\Omega_V 
\neq 
\mathbb{R}^{r \times n}$), the problem becomes difficult in most cases.  
For example, WLRA, robust PCA, NMF, and sparse PCA are all NP hard~\cite{GG10c, 
GV15c, V10, MWA06}. An active direction of research is developing approximation 
algorithms for such problems; see, for example, \cite{CW15} for the norm 
{$\sum_{j=1}^n \|M(:,j)-UV(:,j)\|_2^p$} (for $p=2$, this is PCA), 
\cite{RSW16} for WLRA, and \cite{SWZ16} for the componentwise $\ell_1$-norm. 
\newline

\noindent
{\bf Convexification.}
Under some conditions on the matrix $M$, convexification approaches can lead to 
optimality guarantees. 
When there are no constraints ($\Omega_U = \mathbb{R}^{p \times r}$, $\Omega_V = 
\mathbb{R}^{r \times n}$), 
\eqref{clrma} can be equivalently rewritten as 
\[
\min_{X} \|M-X\| 
\; \text{ such that } \; 
\rank(X) = r . 
\]
From $X$, a solution $(U,V)$ can be obtained by factorizing $X$ (e.g., using 
the SVD). The most widely used convex models are based on minimizing the 
nuclear norm of $X$: 
\begin{equation} \label{clrmaX}
\min_{X} \|M-X\| + \lambda \|X\|_* , 
\end{equation}  
where $\lambda >0$ is a penalty parameter and $\|X\|_* = \sum_{i=1}^{\min(n,p)} 
\sigma_i(X) = \|\sigma(X)\|_1$, $\sigma(X)$ being the vector of singular values 
of $X$. 
This problem can be written as a semidefinite program; see \cite{RFP10} and the 
references therein. 

When the matrix $M$ satisfies some conditions depending on the model (in 
particular, $M$ has to be close to a low-rank matrix), the optimal 
solution to~\eqref{clrmaX} can be guaranteed to recover the solution of the 
original problem; examples include PCA with missing data~\cite{RFP10} and 
robust PCA~\cite{CSPW11,CLM11}. 

As far as we know, these approaches have two drawbacks.
First, if the input matrix $M$ does not satisfy the required conditions, which 
is often the case in practice (e.g., for recommender systems and document 
classification where the input matrix is usually not close to a low-rank 
matrix), 
it is unclear whether the quality of the solution to~\eqref{clrmaX} 
will be satisfactory. 
%
Second, the number of variables is much larger than in~\eqref{clrma}, namely, 
$mn$ 
vs.\@ $r(m+n)$. For large-scale problems, even first-order methods might be too 
costly. 
A possible way to handle the large positive semidefinite matrix is to (re)factor 
it in the SDP as the product of two matrices; this is sometimes referred 
to as the Burer-Monteiro approach~\cite{burer2003nonlinear}.  
In fact, in many cases, any stationary point can be guaranteed to be a global 
minimum~\cite{BVB16, LT16}; see also \cite{lemon2016low} for a survey. 
This is currently an active area of research: trying to identify 
nonconvex problems for which optimal solutions can be guaranteed to be 
computed efficiently (see the end of the next paragraph for other examples). 
\newline

\noindent
{\bf Nonconvex approaches.}
One can tackle~\eqref{clrma} in many ways using standard nonlinear 
optimization schemes. The most straightforward and popular way is to use a 
two-block coordinate descent method (in particular if $\Omega_U$ and $\Omega_V$ 
are convex sets since the subproblems in 
$U$ and $V$ are convex):   
\begin{enumerate} 

\item[0.] Initialize $(U,V)$.

\item $U \leftarrow$ $X$, where $X$  solves exactly or approximately $\min_{X 
\in \Omega_U} \|M-XV\|$.

\item $V \leftarrow$ $Y$, where $Y$ solves exactly or approximately $\argmin_{Y 
\in \Omega_V} \|M-UY\|$.

\end{enumerate}
This simple scheme can be implemented in different ways. 
The subproblems are usually not solved up to high precision; for example, a few 
steps of a (fast) gradient method can be used. 
These methods can in general be guaranteed to converge to a stationary point 
of~\eqref{clrma}~\cite{CR09}. 
More sophisticated schemes include Riemannian optimization 
techniques~\cite{BA11, Bart13}. 
Many methods based on randomization have also been developed recently; see the  
surveys~\cite{mahoney2011randomized, Wood14}.

Alternating and local minimization were shown to lead to optimal solutions under 
 assumptions similar to those needed for convexification-based approaches; 
see, for example, \cite{KOM09, JNS13} for PCA with missing data,  
\cite{agarwal2014learning} for (a variant of) sparse PCA, 
and~\cite{NNSA14} for robust PCA. Recently, \cite{BNS16, GLM16} showed that PCA 
with missing data has no spurious local minima (under appropriate conditions). 
\newline

\noindent
{\bf Outline of the paper.}
In the rest of this paper, we focus on a particular CLRMA problem, namely, 
nonnegative  
matrix factorization (NMF), with 
$\|\cdot \|=\|\cdot \|_F^2$,  
$\Omega_U = \mathbb{R}^{p \times r}_+$, and 
$\Omega_V = \mathbb{R}^{r \times n}_+$. 
As opposed to other CLRMA variants (such as robust PCA, sparse PCA, and PCA 
with 
missing data), as far as we know, no useful convexification approach exists.

The goal of this paper is not to provide an exhaustive survey but rather to 
provide a brief introduction, focusing only on several aspects of NMF (obviously 
biased toward our own interests). In particular, we address 
the application of NMF for hyperspectral imaging,  
the geometric interpretation of NMF, 
complexity issues, 
algorithms, and 
the nonnegative rank and its link with extended formulations of polyhedra.

\section{Nonnegative Matrix Factorization} 

The standard NMF problem can be formulated as follows 
\begin{equation} \label{nmf}
\min \limits_{U \in \mathbb{R}^{p \times r}, V\in \mathbb{R}^{r \times n}} 
\|M-UV\|_F^2 \,\,
\text{ such that} \,\,
U, V \geq 0.
\end{equation}
As mentioned in the introduction, these nonnegativity constraints allow 
interpreting the basis elements in the same way as the data (e.g., as image, or 
vector of word counts) while the nonnegativity of $V$ allows interpreting the 
weights as activation coefficients. 
We describe in detail in the next section a particular 
application, namely,
blind hyperspectral unmixing, where the nonnegativity of $U$ and $V$ has a 
physical interpretation. 

The nonnegativity constraints also naturally lead to sparse factors. In fact, 
the first-order optimality conditions of a problem of the type  
$\min_{x \geq 0} f(x)$ are $x_i \geq 0$, $\nabla_i f(x) \geq 0$ and $\nabla_i 
f(x) x_i = 0$ for all $i$. Hence stationary points 
of~\eqref{nmf} are expected to have zero entries. 
This property of NMF enhances its interpretability and provides a better 
compression compared with unconstrained variants.  

We refer to the problem of finding an exact factorization, that is, finding 
$U \geq 0$ and $V \geq 0$ such that $M=UV$, as ``exact NMF.'' 
The minimum $r$ such that an exact NMF exists is the nonnegative rank of $M$, 
denoted $\rank_+(M)$. 
We have that $\rank(M) \leq \rank_+(M) \leq \min(m,n)$ (since $M = MI = IM$,  
where $I$ is the identify matrix). 

NMF has been used successfully in many applications; 
see, for example, \cite{CZA09, G14}  and the references therein. 
In the next section we focus on one particular application, namely, blind 
hyperspectral unmixing.

\section{Hyperspectral Imaging} \label{blindHU}
\begin{figure}[ht!] 
\begin{center}
\includegraphics[width=8.0cm]{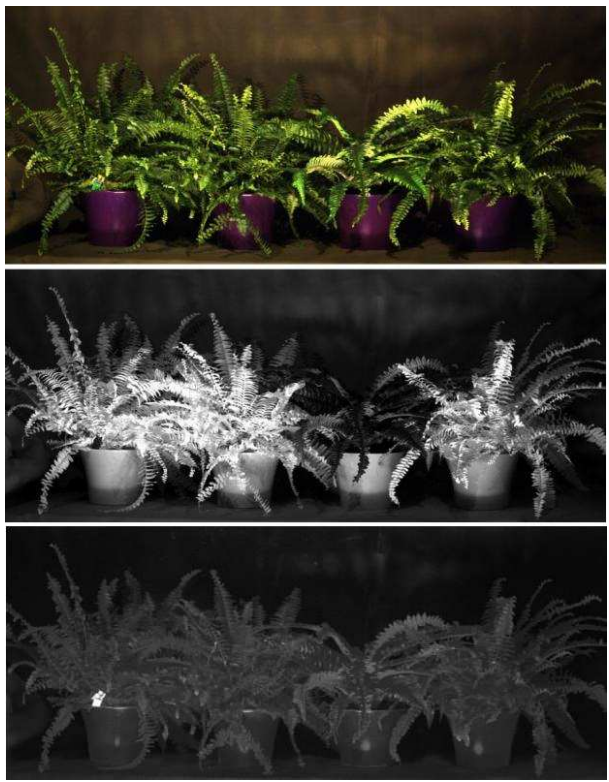}   
\caption{From top to bottom: 
(1)~RGB image of four plants: can you identify the artificial one? 
(2)~Grayscale image at a wavelength that is blind to the naked eye (namely,  
770 nm, infrared) and allows identifying the artificial plant 
(plants have a high reflectance at infrared wavelengths, as opposed to the 
artificial material). 
(3)~Analysis of the image allows finding a small target, a LEGO figure within 
the plants.  
Source: sciencenordic.com, Photo: Torbj{\o}rn Skauli, FFI. \label{plants}} 
\end{center}
\vspace*{-.5cm}
\end{figure}

A grayscale image is an image in which the value of each pixel is a single 
sample. An RGB image has three channels (red, green, and blue) and allows a 
color image to be reconstructed as it is perceived by an human eye. 
A hyperspectral image is an image for which usually each pixel has between 100 
and 200 channels, 
corresponding to the reflectance (fraction of light reflected by that pixel) at 
different wavelengths.  
The wavelengths measured in a hyperspectral image depend on the camera used and 
are usually chosen depending on the application at hand. 
The advantage of hyperspectral images is that they contain much more 
information, some of it blind to the human eye, that allows one to identify and 
characterize the materials present in a scene much more precisely; see 
Figure~\ref{plants} for an illustration. 
Its numerous applications include agriculture, 
	eye care, 
	food processing, 
	mineralogy, 
	surveillance, 
	physics, 
	astronomy, 
	chemical imaging, and 
	environmental science; 
see, for example, 
\url{https://en.wikipedia.org/wiki/Hyperspectral_imaging} or 
\url{http://sciencenordic.com/lengthy-can-do-list-colour-camera}.

Assume a scene is being imaged by a hyperspectral imager using $p$ wavelengths 
(that is, $p$ channels) and $n$ pixels. 
Let us construct the matrix $M \in \mathbb{R}^{p \times n}_+$ such that 
$M(i,j)$ is the reflectance of the $j$th pixel at the $i$th wavelength. 
Each column of $M$ therefore corresponds to the so-called spectral signature of 
a pixel, while each row corresponds to a vectorized image at a given 
wavelength. 
Given such an image, an important goal in practice is to (1) identify the 
constitutive materials present in the image, called endmembers (e.g., grass, 
trees, road surfaces, roof tops) and (2) classify the pixels accordingly, that 
is, identify which pixels contain which materials and in which quantity. In 
fact, the resolution of most hyperspectral images is low, and hence most 
pixels will contain several materials. 
If a library or dictionary of spectral signatures of materials present in the 
image 
is not available, this problem is referred to as blind hyperspectral unmixing 
(blind HU): the goal is to identify the endmembers and quantify the abundances 
of the endmembers in each pixel. 

The simplest and most popular model is the linear mixing model (LMM). It assumes 
that the spectral signature of a pixel equals the weighted linear combination of 
the spectral signatures of the endmembers it contains, where the weight is  
given by the abundances. 
Physically, the reflectance of a pixel will be proportional to the materials it 
contains: for example, if a pixel contains 30\% of aluminum and 70\% of copper, 
its spectral signature will be equal to 0.3 times the spectral signature of the 
aluminum plus 0.7 times the spectral signature of the copper. In practice, this 
model is only approximate because of imperfect conditions (measurement noise, 
light reflecting off several times before being measured, atmospheric 
distortion, 
etc.). We refer the reader to~\cite{BP12, Ma14} for recent surveys on (blind) HU 
techniques and to~\cite{smith2006} for an introduction to hyperspectral imaging.

\begin{figure}[ht!] 
\begin{center}
\includegraphics[width=12.0cm]{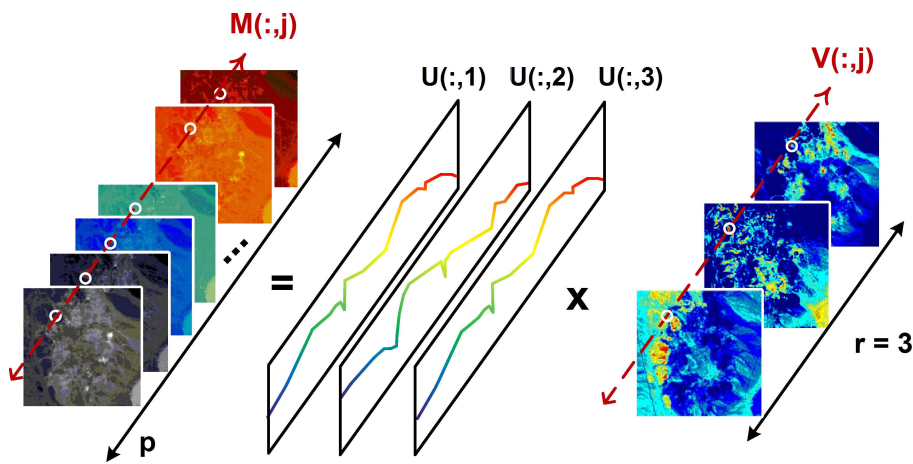} 
\caption{Illustration of the decomposition of a hyperspectral image with three 
endmembers~\cite{Ma14}. 
On the left, the hyperspectral image $M$; in the middle, the spectral 
signatures 
of the three endmembers as the columns of matrix $U$; on the right, the 
abundances of each material in each pixel (referred to as the abundance maps). 
\label{bhu}} 
\end{center}
\vspace*{-.5cm}
\end{figure}

If we use the LMM and assume that the image contains $r$ endmembers whose 
spectral 
signatures are given by the columns of the matrix $U \in \mathbb{R}^{m \times 
r}_+$, we have for all $j$ 
\[
M(:,j) = \sum_{k=1}^r v_{kj} U(:,k) = U V(:,j), 
\]
where $v_{kj} \geq 0$ is the abundance of the $k$th endmember in the $j$th 
pixel. Therefore, blind HU boils down to the NMF of matrix $M$; see 
Figure~\ref{bhu} for an illustration. 

Using a standard NMF algorithm, that is, an algorithm that tries to 
solve~\eqref{nmf}, will in general not lead to the sought decomposition. The 
reason is that the solution of NMF is highly nonunique, as discussed later. 
In practice,
a meaningful solution
is achieved usually by  
using additional constraints/penalty terms, including: the sum-to-one 
constraints on the 
abundances ($\sum_{k=1}^r v_{kj} = 1 \, \forall j$), sparsity of $V$ (because 
most pixels contain only a few endmembers), 
piecewise smoothness of the columns of $U$ (since they correspond to spectral 
signatures), and
spatial coherence of the rows of $V$ (because neighboring pixels are more likely 
to contain the same endmembers).  
Numerous constrained variants of NMF exist that we do not discuss here; 
see, for example, \cite{CZA09, G14} and the references therein.

\section{Geometry and Uniqueness} \label{geocou}

NMF has a nice geometric interpretation, which is crucial to consider in order 
to 
understand the nonuniqueness of the solutions. 
As discussed subsequently, it also allows one to develop efficient algorithms
and is closely related to the extended formulations of polyhedra.
  
Let us consider the exact case, that is, $M = UV$. 
Without loss of generality, (i) the zero columns of $M$ and $U$ can be removed, 
and (ii) the columns of $M$ and $U$ can be normalized so that the entries of 
each column sum to one: 
\[
M D_M^{-1} = U D_U^{-1} D_U V D_M^{-1}, 
\] 
where $D_M$ and $D_U$ are diagonal matrices with $D_M(j,j) = \|M(:,j)\|_1$ and 
$D_U(j,j) = \|U(:,j)\|_1$, respectively. 
Since we have 
$M(:,j) = \sum_{k=1}^r U(:,k) V(k,j) = U V(:,j)$,  
this normalization implies that the columns of $V$ also have their entries 
summing to one, that is, $\|V(:,j)\|_1 = 1$ for all $j$. 
Thus that, after normalization, the columns of $M$ belong to the convex 
hull of the columns of $U$:
\[
M(:,j)  \, \in \, 
\conv(U) \, \subseteq \, \Delta^p = \{ x \in \mathbb{R}^{p} | x \geq 0, 
\|x\|_1 = 1\} \quad \forall j,
\] 
where $\conv(U)=\{ Ux | x \geq 0, \|x\|_1 = 1\}$. Therefore, the exact NMF 
problem is equivalent to finding a polytope, $\conv(U)$, nested between two 
given polytopes,  $\conv(M)$ and the unit simplex $\Delta^p$. 
The dimension of the inner polytope, $\conv(M)$, is $\rank(M)-1$, while the 
dimension of the outer polytope, $\Delta^p$, is $p-1$. The dimension of the 
nested polytope $\conv(U)$ is not known in advance. 
When the three polytopes (inner, nested, and outer) have the same dimension, 
this 
problem is well known in computational geometry and is referred to as the 
nested 
polytope problem (NPP)~\cite{DJ90}. 

If $\rank(M) = \rank(U)$, the column spaces of $M$ and $U$ must coincide, and 
the 
outer polytope can be restricted to $\Delta^p \cap \col(M)$,  
 in which case the inner, nested, and outer polytopes have the same dimension. 
If we impose explicitly this additional constraint ($\rank(M) = \rank(U)$) on 
the 
exact NMF problem, 
we can prove that NPP and this restricted variant of exact NMF are 
equivalent, that is, they can be reduced 
to one another~\cite{GG12nr,CK16}.

To illustrate, we present a simple example with nested hexagons; this 
is similar 
to the example presented in~\cite{MSV03}. Let $a > 1$, and let $M_a$ be the 
matrix 
\begin{equation} \label{hexa}
\frac{1}{a} \left( \begin{array}{cccccc}
1    & a    & 2a-1 & 2a-1 & a & 1 \\
1    & 1    & a    & 2a-1 & 2a-1  & a  \\
a    & 1    & 1    & a & 2a-1  & 2a-1  \\
2a-1 & a    & 1    & 1 & a & 2a-1  \\
2a-1 & 2a-1 & a    & 1 & 1 & a \\
a    & 2a-1 & 2a-1 & a & 1 & 1 \\
 \end{array}  \right) . 
\end{equation}
The restricted exact NMF problem for $M_a$ involves two nested hexagons 
(recall that 
we restrict the polytopes to be in 
the intersection between the column space of $M_a$ and $\Delta^p$, which has 
dimension 2 since $\rank(M_a) = 3$). 
Each facet of the outer polytope corresponds to a facet of the 
nonnegative orthant, 
that is, to a nonnegativity constraint. 
The inner hexagon is smaller than the outer one with a ratio of $\frac{a-1}{a}$. 

\begin{figure}[t!] 
\begin{center}
\begin{tabular}{cc} 
\includegraphics[width=8.35cm]{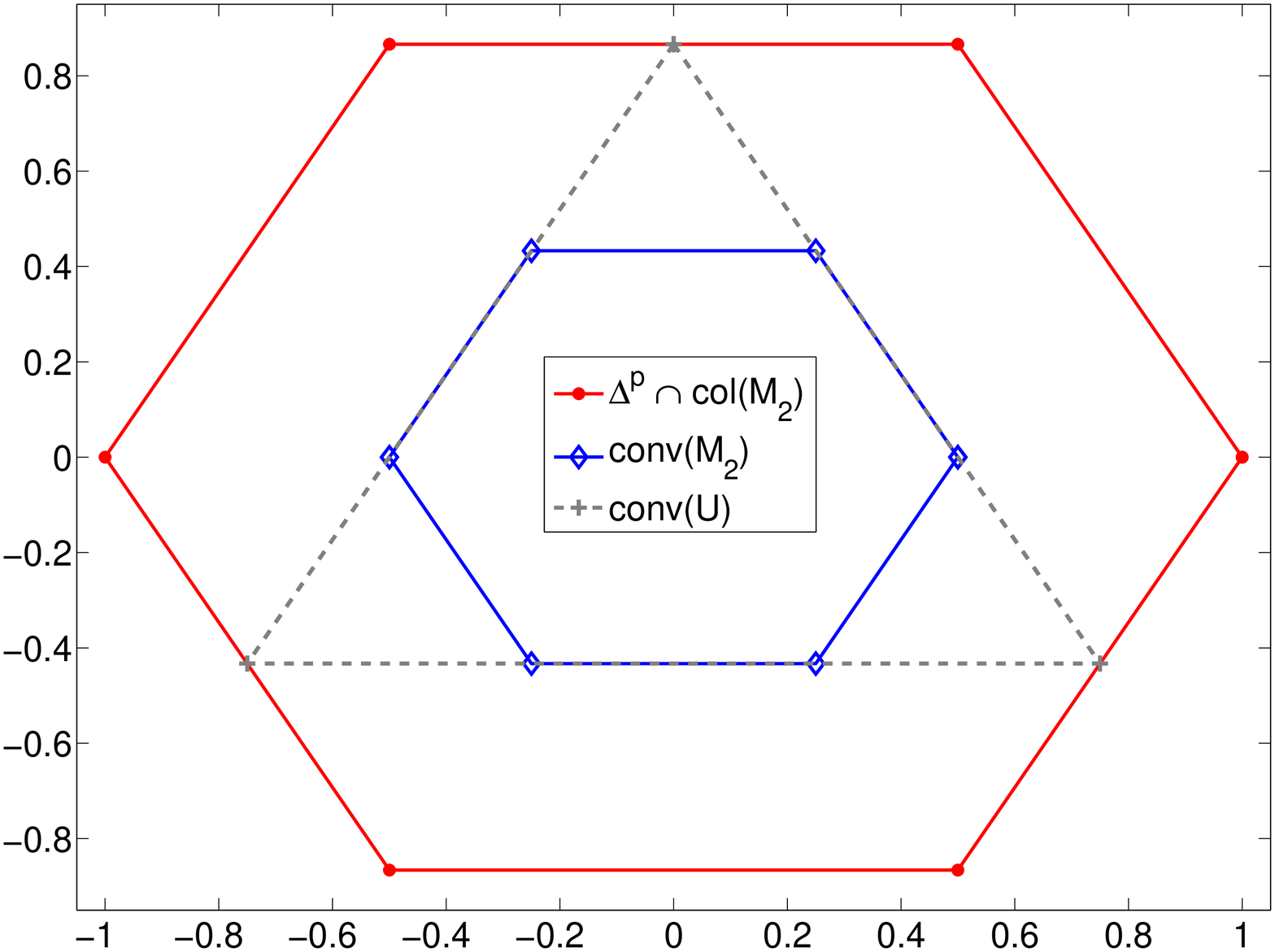} & 
\includegraphics[width=8.35cm]{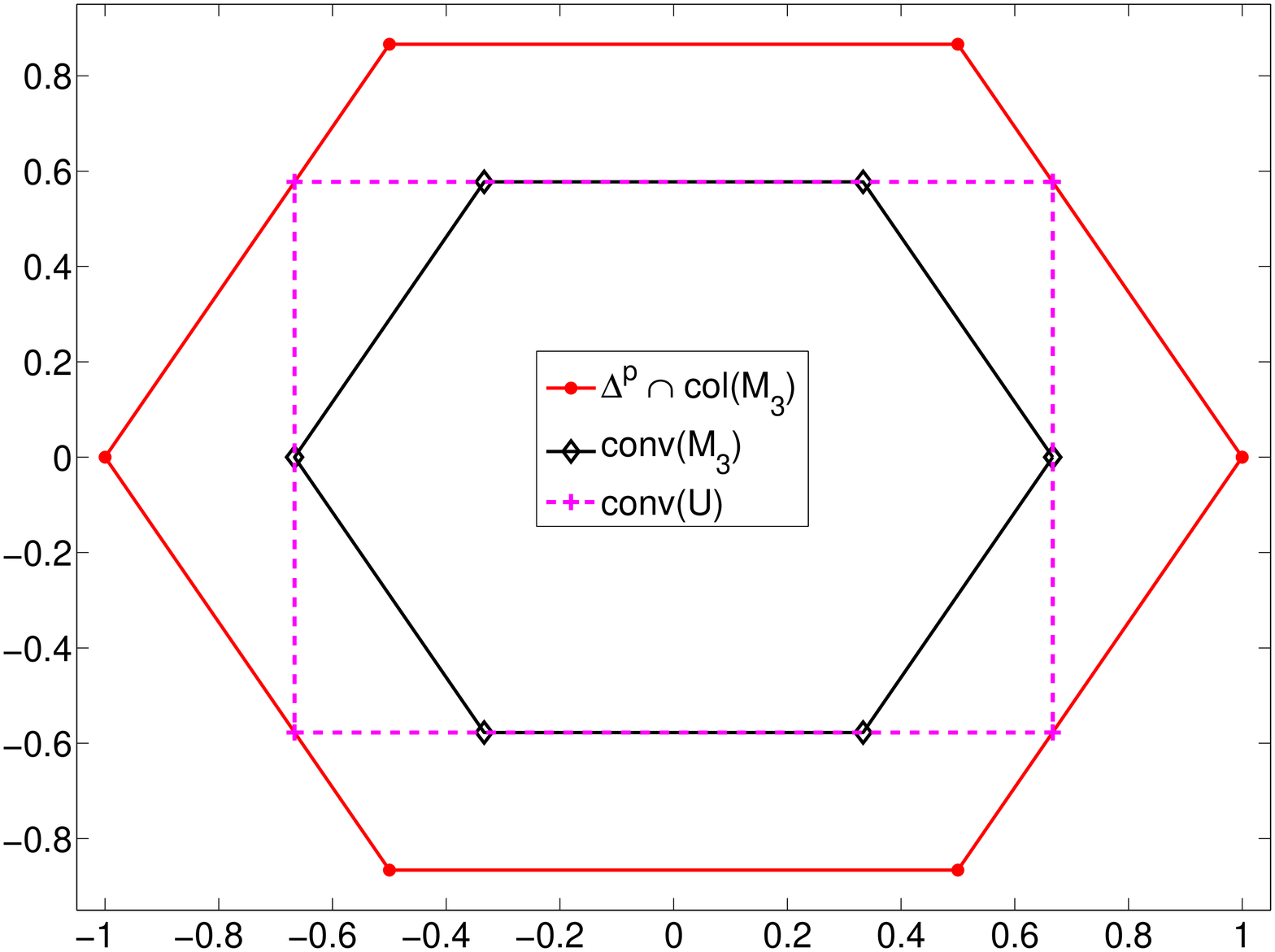} 
\end{tabular}
\caption{NPP problem corresponding to the exact NMF of the matrix 
from~\eqref{hexa}, restricted to the column space of $M$: (top) the case 
$a = 2$; (bottom) $a = 3$. \label{nesthex}} 
\end{center}
\vspace*{-.5cm}
\end{figure}

For $a=2$, the inner hexagon is twice as small as the outer one, and we can fit 
a triangle between the two so that $\rank_+(M_a) = 3$; 
see Figure~\ref{nesthex} (top). For any $a > 2$, $\rank_+(M_a) \geq 4$ because 
no triangle can fit between the two hexagons.  
For $a=3$, the inner hexagon is 2/3 smaller than the outer one, and we can fit 
a 
rectangle between the two and $\rank_+(M_a) = 4$; 
see Figure~\ref{nesthex} (bottom). This implies that $\rank_+(M_a) = 4$ for all 
$2 < a \leq 3$. 

\begin{figure}[htb!] 
\begin{center}
\includegraphics[width=13.0cm]{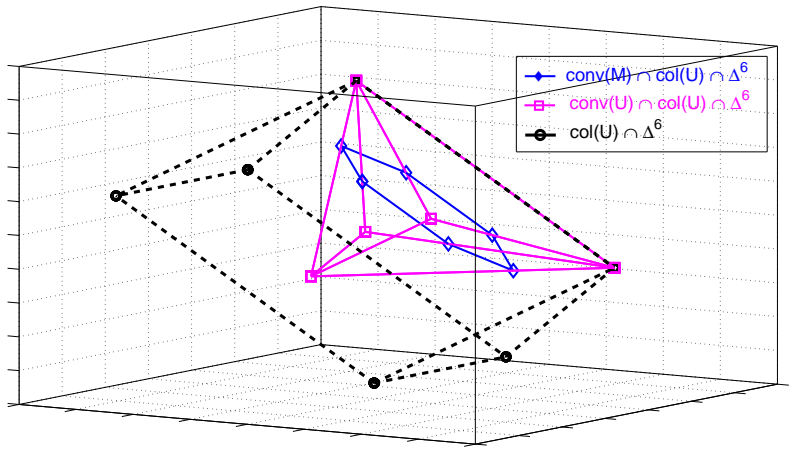} 
\caption{NPP solution corresponding to the exact NMF of the matrix 
from~\eqref{neshexa}, restricted to the column space of $U$. 
It corresponds to the matrix $M_a$ from~\eqref{hexa} when $a \rightarrow 
\infty$. \label{nesthex3d}} 
\end{center}
\vspace*{-.5cm}
\end{figure}

For any $a > 3$, $\rank_+(M_a) = 5$. Surprisingly, the nonnegative rank of $M_a$ 
is always 
no more 
than 5 (even when $a$ tends to infinity, in which case the 
inner and outer hexagons coincide) because there exists a three-dimensional 
polytope 
within $\Delta^6$ with 5 vertices that contains the outer polytope; see 
Figure~\ref{nesthex3d}, which corresponds to the factorization  
\begin{equation} \label{neshexa}
M =  \lim_{a \rightarrow +\infty} 
M_a = \left( \begin{array}{cccccc}
         0  &   1 &    2   &  2   &  1  &        0 \\
         0      &   0 &   1   &  2   &  2   &  1  \\
    1   &  0     &     0  &  1  &   2 &    2  \\
    2    & 1   &       0    &     0 &   1   &  2  \\
    2   &  2   &  1 &    0  &   0   &  1  \\
    1   &  2 &    2  &   1  &   0   &       0 \\
		\end{array} \right) \\
= UV  =  
		\left( \begin{array}{ccccc}
	   1  &   0  &   0 &    1  &   0 \\
     2  &   0  &   0 &    0  &   1 \\
     1   &  0   &  1  &   0  &   0 \\
     0   &  1   &  1  &   0   &  0 \\
     0   &  2   &  0  &   0   &  1 \\
     0  &   1  &   0  &   1  &   0  \\
		\end{array} \right) 
		\left( \begin{array}{cccccc}
		 0  &   0  &   0   &  1 &   1   &  0\\ 
     1  &   1  &   0   &  0 &    0  &   0 \\
     1  &   0  &   0   &  0 &   1   &  2 \\
     0   &  1   &  2    & 1  &   0   &  0 \\
     0   &  0  &   1    & 0   &  0  &   1 \\ 
				\end{array} \right), 
\end{equation}
where $\rank(U) = 4$, and hence $\conv(U)$ has dimension 3. 

This example illustrates other interesting properties of NMF: 
\begin{itemize} 

\item NMF does not in general have a unique solution (up to scaling and 
permutation of the rank-one factors). 
For example, for $a=2$ (Figure~\ref{nesthex}, top), four triangles 
 can be fit between the two polytopes 
(the one shown on the figure, its rotation by 60 degrees, and two triangles 
whose vertices are three nonadjacent vertices 
of the outer hexagon).
For $1<a<2$, this would be even worse since there would be an infinite number of 
solutions. 
For this reason, practitioners often add additional constraints to the NMF model 
to try to identify the most meaningful solution to their problem (such as 
sparsity, as discussed earlier);
see, for example, \cite{LCP08, G12, HSS14} 
for more details on the uniqueness of NMF. 

\item The nonnegative rank can increase only in the neighborhood of a given 
matrix; that is, the nonnegative rank is upper 
semicontinuous~\cite[Th.3.1]{BCR11}: ``If $P$ is a nonnegative matrix, without zero columns and with $\rank_+(P) = r$,
then there exists a ball $B(P, \epsilon)$ centered at $P$ and of radius 
$\epsilon > 0$ such that $\rank_+(N) \geq r$ for all $N \in B(P, e)$.'' 

\end{itemize}

\section{Complexity} \label{complex}

Given a nonnegative matrix $M$, checking whether $\rank(M) = \rank_+(M) = r$ is 
NP hard: unless $P=NP$, there is no polynomial time algorithm in $m$, $n$ and 
$r$ for this problem~\cite{V10}. 
If $r$ is fixed, however, there is a polynomial time algorithm running in 
$O\big( 
(pn)^{r^2} \big)$~\cite{AGKM11, Moit13}.
 The argument is based on quantifier elimination theory 
(in particular the fact that checking whether a system of $\ell$ equations in 
$n$ variables up to degree $d$ can be solved in time polynomial in $\ell$ and 
$d$).  
Unfortunately, as far as we know, this cannot be used in practice, even for 
small matrices (e.g., checking whether a 4-by-4 matrix has nonnegative 
rank 3 seems already impractical with current solvers). 
Developing an effective code for exact NMF for small matrices is an 
important direction for further research. Note that we have developed a code 
based on heuristics that allows solving exact NMF for matrices up to a few 
dozen rows and columns (although our code comes with no 
guarantee)~\cite{vandaele2016heuristics}.  

More recently, Shitov~\cite{shitov2016} and independently Chistikov et 
al.~\cite{CKMS16} answered an important open problem showing that the 
nonnegative rank over the reals might be different from the nonnegative rank 
over the rationals, implying that the nonnegative rank computation is not in NP 
since the size of the output is not bounded by the size of the input.

\section{Algorithms} \label{algoNMF}

In this section, we briefly describe the two main classes of NMF algorithms. As 
mentioned in the introduction, there does not exist, to the best of our 
knowledge, a successful 
convexification approach for NMF, as opposed to other low-rank models. Note, 
however, that there does exist a convexification approach to 
compute lower bounds for the nonnegative rank~\cite{Fawzi2015}. 
An explanation is that we cannot work directly with the low-rank 
approximation $X=UV$ and use the nuclear norm of $X$, because even if we were 
given the best nonnegative approximation $X$ of nonnegative rank $r$ for $M$, 
in general recovering the exact NMF $(U,V)$ of $X$ would be difficult.
Writing directly a convexification in variables $(U,V)$ seems difficult (for 
rank higher than one\footnote{Note that the rank-one NMF problem is equivalent 
to the rank-one unconstrained problem since for any rank-one solution $uv^T$, 
one can easily check that $|u\|v|^T$ is a solution with lower objective 
function 
value. 
This also follows from the Perron-Frobenius and Eckart-Young theorems.}) 
because of the symmetry of the problem (permuting columns of $U$ and rows of $V$ 
accordingly provides an equivalent solution). 
Breaking this symmetry seems nontrivial; see~\cite[pp.\ 146-148]{NG11} for a 
discussion and a tentative SDP formulation.  
This is an interesting direction for further research.

\subsection{Standard nonlinear optimization schemes}  

As for CLRMA problems, most NMF algorithms use a two-block coordinate descent 
scheme: 
\begin{enumerate} 

\item[0.] Initialize $(U,V) \geq 0$.

\item $U \leftarrow$ $X$, where $X$  solves exactly or approximately $\min_{X 
\geq 0} \|M-XV\|_F$ .

\item $V \leftarrow$ $Y$, where $Y$ solves exactly or approximately $\argmin_{Y 
\geq 0} \|M-UY\|_F$ .

\end{enumerate} 
Note that the subproblems to be solved are so-called nonnegative least squares 
(NNLS). Because NMF is NP hard, these algorithms can only guarantee convergence 
(usually to a first-order stationary point). 

The most well-known algorithm for NMF is the multiplicative updates, 
namely, 
\[
U \leftarrow U \circ \frac{[ M V^T ]}{[ UVV^T ]}, \quad V \leftarrow V \circ 
\frac{[ U^T M ]}{[ U^T UV ]} , 
\]
where $\circ$ (resp.\@ $\frac{[ \, ]}{[ \, ]}$) is the componentwise product 
(resp. division) between two matrices. 
It is extremely popular because of its simplicity and because it was proposed  
in the paper of Lee and Seung~\cite{LS99} that launched the research on NMF.  
However, it converges slowly; it cannot modify zero entries; and it is not 
guaranteed to converge to a stationary point. Note that it can be interpreted as 
a rescaled gradient descent; see, for example,~\cite{G14}.  

Methods that try to solve the subproblems exactly are referred to as alternating 
nonnegative least squares; among these, active set methods seem to be 
the most 
efficient, and dedicated codes have been implemented by Haesun Park and 
collaborators; see~\cite{kim2014algorithms} and the references therein.  

In practice, a method that seems to work extremely well is to apply a few steps 
of coordinate descent on the NNLS subproblems: the subblocks are the columns of 
$U$ and the rows of $V$ \cite{CZA07, GG12}---the reason is that the 
subproblems 
can be solved in closed form. In fact, the optimal $k$th column of $U$ (all 
other variables being fixed) is given by 
\[ \displaystyle
\arg\min \limits_{U(:,k) \geq 0} \| R_k - U(:,k) V(k,:)\|_F^2 
= 
\max \left( 0 , \frac{ R_k V(k,:)^T  }{\|V(k,:)\|_2^2} \right ) , 
\]
for $R_k = M - \sum_{j \neq k} U(:,j) V(j,:)$, and similarly by symmetry for 
the $k$th row of $V$. 

Many other approaches can be applied to the NNLS subproblems (e.g., 
projected gradient method~\cite{lin2007projected}, 
fast/accelerated gradient method (Nesterov's method)~\cite{guan2012nenmf}, 
and Newton-like method~\cite{CZA06}).

\subsection{Separable NMF}  

Although they usually provide satisfactory results in practice, the methods 
described in the preceding section do not come with any guarantee. 
In their paper on the complexity of NMF, Arora et al.~\cite{AGKM11} also 
identify a subclass of matrices for which the NMF problem is much easier. These 
are the so-called separable matrices defined as follows. 
\begin{definition} A matrix $M$ is separable if there exists a subset 
$\mathcal{K}$ of $r$ of its columns with $r = \rank_+(M)$ and a nonnegative 
matrix $V$ such that 
$M = M(:,\mathcal{K}) V$. 
\end{definition}
This requires each column of the basis matrix $U$ in an NMF decomposition to be 
present in the input matrix $M$. 
Equivalently, this requires the matrix $V$ in an NMF decomposition to contain 
the identity matrix as a submatrix. The separable NMF problem is the problem to 
identify the subset $\mathcal{K}$ (in the noisy case, this subset should be such 
that 
$\min_{V \geq 0} \|M - M(:,\mathcal{K}) V\|$ is minimized). 

Although this condition is strong, it makes sense in several 
applications, for example the following.
\begin{itemize} 

\item Document classification: for each topic, there is a ``pure'' word used 
only by that topic (an ``anchor'' word)~\cite{Ar13}.  

\item Time-resolved Raman spectra analysis: each substance has a peak in its 
spectrum while the other spectra are (close to) zero~\cite{LHKL16}.  

\item Blind hyperspectral unmixing: for each endmember, there exists a pixel 
that contains only that endmember. This is the so-called pure-pixel assumption 
that has been used since the 1990s in that community. 

\end{itemize}  
Other applications include video summarization~\cite{ESV12} and 
foreground-background separation~\cite{KS15}. 

Geometrically, in the exact case and after normalization of the columns of $X$ 
and $U$, 
the separability assumption is equivalent to having $\conv(U) = \conv(M)$. 
Therefore, the so-called separable NMF problem reduces to identify the vertices 
of the convex hull of the columns of $M$. This is a relatively easy geometric 
problem. 
It becomes tricky when noise is added to the separable matrix, and many recent 
works have tried to quantify the level of noise that one can tolerate and  
still be able to recover the vertices, up to some error.

\subsubsection{Geometric algorithms} \label{geoalgo}

Most algorithms for separable NMF are based on the geometric interpretation,  
many being developed within the blind HU community (sometimes referred to as  
pure-pixel search algorithms). 
Only recently, however, has robustness to noise of these algorithms been 
analyzed.  

 One of the simplest algorithm, often referred to as the successive projection  
algorithm, is closely related to the modified Gram-Schmidt algorithm with 
column pivoting and has been discovered several times~\cite{MC01, RC03, CM11}; 
see the discussion in~\cite{Ma14}. 
Over a polytope, a strongly convex function (such as the $\ell_2$ norm) is 
always  maximized at a vertex: this can be used to identify a vertex, that is, a 
column of $U$ (recall that we assume that the columns of $M$ are normalized so 
that $\conv(U) = \conv(M)$ under the separability assumption). Once a column of 
$U$ has been identified, one can project all columns of $M$ onto the orthogonal 
complement of that column (so that this particular column projects onto 0): this 
amounts to applying a linear transformation to the polytope. 
If $U$ is full rank (meaning the polytope is a simplex, which is the case 
usually in practice), 
then the other vertices do not project onto 0, and one can use these two steps  
recursively. This approach is a greedy method to identify a subset of the 
columns with maximum volume~\cite{AM09, CM11}. 
This algorithm was proved to be robust to noise~\cite{GV14} and can be made  
more robust to noise by using strategies such as
\begin{itemize}

\item applying dimensionality reduction, such as PCA, to the columns of $M$ in 
order to filter the noise~\cite{ND05}; 

\item using a preconditioning based on minimum-volume ellipsoid~\cite{GV15, 
Mizu2016}; 

\item going over the identified  
vertices (once $r$ vertices have been identified) to check whether they still 
maximize the strongly convex function once 
projected onto the orthogonal complement of the other vertices (otherwise, they 
are replaced, increasing the volume of the identified 
vertices)~\cite{Ar13}; 

\item taking into account the nonnegativity constraints in the projection  
step~\cite{Nic14}. 

\end{itemize}

We refer the reader to~\cite{BP12, Ma14} for surveys on these approaches. 
Most geometric approaches for separable NMF are computationally cheap. 
Usually, however, they are sensitive to outliers.

\subsubsection{Convex models}

If $M$ is separable, there exist an index set $\mathcal{K}$ of size $r$ 
and a nonnegative matrix $V$ such that $M = M(:,\mathcal{K})V$. Equivalently,  
there exists an $n$-by-$n$ nonnegative matrix $X$ with $r$ nonzero rows such 
that $M = MX$ with $X(\mathcal{K},:) = V$. Solving separable NMF can therefore 
be formulated as 
\[
\min_{X \geq 0} \|X\|_{\text{row},0} \quad \text{ such that } M = MX, 
\] 
where $\|X\|_{\text{row},0}$ counts the number of nonzero rows of $X$. A 
standard  convexification approach is to use the $\ell_1$ norm, 
replacing $\|X\|_{\text{row},0}$ with $\sum_{i=1}^n \| X(i,:) \|_{k}$ for some  
$k$; for example,  
\cite{ESV12} uses $k=\infty$ and \cite{EMO12} uses $k=2$. 

If the columns of $M$ are normalized, the entries of $V$ are bounded above by  
one (since the columns of $U$ are vertices), and 
another formulation for separable NMF is obtained:
\[
\min_{X \geq 0}  \|\diag(X)\|_{0} \; \text{ such that } \; M = MX \text{ and } 
X(i,j) \leq X(i,i) \leq 1 \, i,j.
\]
Because  
on each row the diagonal entry has to be the largest and because the goal is to 
minimize the number of nonzero entries of the diagonal of $X$, the optimal 
solution will contain $r$ nonzero diagonal entries and hence $r$ nonzero rows. 
(Note that requiring the diagonal entries of $X$ to be binary would allow one 
to 
model this problem exactly by using mixed-integer linear programming.) 
Using the $\ell_1$ norm, we get another convex model (proposed 
in~\cite{BRRT12} 
and improved in~\cite{GL14}): 
 \[ 
\min_{X \geq 0}   \tr(X)  
\text{ such that} M = MX \; \text{ and } \; 
X(i,j) \leq X(i,i) 
\leq 
1 \, \forall i,j, 
\] 
where $\tr(X)$ 
is equal to $\|\diag(X)\|_{1}$ since $X$ is nonnegative. 
In practice, when noise is present, the equality term $M = MX$ is replaced with 
$\|M - MX\| \leq \epsilon$ for some appropriate norm 
(typically the $\ell_1$, $\ell_2$, or Frobenius norm)  
or is added in the objective function as a penalty. 

The two models presented above turn out to be essentially 
equivalent~\cite{GL16}. 
The main drawback is the computational cost, since these models have $n^2$ 
variables. 
For example, in hyperspectral imaging, $n$ is the number of pixels and is 
typically on  
the order of millions; hence, solving these problems is challenging (if 
not impractical). 
A natural approach is therefore to first select a subset of good candidates 
among the columns of $M$ (e.g., using geometric algorithms) and then  
optimize only over this subset of the rows of $X$~\cite{EMO12, GL16}.  
The main advantage of this approach is that the resulting models are 
provably the most robust for 
separable NMF~\cite{GL14}. 
Intuitively, the reason is not only that the model focuses in 
identifying, 
for 
example, a subset of columns with large volume but also that it requires all 
the data 
points to be well approximated with the selected vertices (since $\|M-MX\|$ 
should be small). 
For this reason, they are also much less sensitive to outliers than are most 
geometric approaches.

\section{Nonnegative Rank and Extended Formulations}  \label{extform}

We now describe the link between extended formulations of polyhedra and NMF. 
This is closely related to the geometric interpretation of 
NMF described earlier. Let $\mathcal{P}$ be a polytope 
\[
\mathcal{P} = \{ x \in \mathbb{R}^k \ | \ b_i - A(i,:)x \geq 0  \text{ for } 1 
\leq i \leq p \}, 
\]
and let ($w_1, \, \cdots, \, w_n$) be its vertices. Let 
$S_\mathcal{P}$ be the 
$p$-by-$n$ slack matrix  of $\mathcal{P}$ defined as follows: 
\[
S_\mathcal{P}(i,j) = b_i - A(i,:)w_j 
\qquad 1 \leq i \leq p, \, 1 \leq j \leq n. 
\]
An extended formulation of $\mathcal{P}$ is a higher-dimensional polyhedron $Q 
\subseteq \mathbb{R}^{k+p}$ that (linearly) projects onto $P$. 
The minimum number of facets (that is, inequalities) of such a polytope is 
called the extension complexity, xp($\mathcal{P}$), of $\mathcal{P}$. 
\begin{theorem} (Yannakakis, \cite{Y91}). \label{th1} 
Let $S_\mathcal{P}$ be the slack matrix of the polytope $\mathcal{P}$. Then, 
$
 \rank_+(S_\mathcal{P}) 
 = 
\textrm{xp}(\mathcal{P}). 
$ 
\end{theorem} 
Let us just show that $\text{xp}(\mathcal{P}) \leq 
\rank_+(S_\mathcal{P})$, because it is elegant and straightforward. 
Given $\mathcal{P} = \{ x \in \mathbb{R}^k \ | \ b - A x \geq 0  \}$, 
any exact NMF of $S_\mathcal{P} = UV$ with $U \geq 0$ and $V \geq 0$ 
provides an explicit extended formulation (with some redundant equalities) 
of $\mathcal{P}$: 
\[  
Q = \{ (x,y) \ |\ b - Ax = U y \text{ and } y \ge 0 \} . 
\]
In fact, let us show that $Q_x = \{ x | \exists y \text{ s.t. } (x,y) \in Q\} 
= \mathcal{P}$. 
We have $Q_x \subseteq \mathcal{P}$ since $U \geq 0$ and $y \geq 0$; hence $b - 
Ax = U y \geq 0$ for all $(x,y) \in Q$.  
We have $\mathcal{P} \subseteq Q_x$ because all vertices of $\mathcal{P}$ belong 
to $Q_x$: by construction, 
$(w_j,V(:,j)) \in Q$ since $S_\mathcal{P}(:,j) = b - A w_j = U V(:,j)$ and 
$V(:,j) \geq 0$. 
\newline

\noindent
{\bf Example.} The extension complexity of the regular 
$n$-polygons is $O(\log_2(n))$~\cite{FRT12}. 
This result can be used to approximate a second-order cone program with a 
linear program~\cite{BTN01}. 
In particular, we have seen 
that the extension 
complexity of the regular hexagon is 5; see Equation~\eqref{neshexa} and 
Figure~\ref{nesthex3d}. 
\newline

\noindent
{\bf Recent results.}
Several recent important results for understanding the limits of linear 
programming for solving combinatorial problems
are based on Theorem~\ref{th1} and on constructing lower bounds for the 
nonnegative rank, usually based on the sparsity pattern of the slack 
matrix~\cite{fiorini2013combinatorial}; see~\cite{K11} for a survey. 
In particular, Rothvo{\ss} showed recently that the prefect matching problem 
cannot be written with polynomially many 
constraints~\cite{rothvoss2014matching}. 

 These ideas can be generalized in two ways:  
\begin{itemize} 

\item To characterize the size of approximate extended formulations (for a 
given precision)~\cite{BFPS15}. 
 
\item To any convex cone~\cite{gouveia2013lifts}, which leads to other 
CLRMA problems. 
For example, for the cone of positive semidefinite (PSD) matrices, the rows 
of $U$ and the columns of $V$ are required to be vectorized PSD matrices. The 
smallest PSD extension of a given set (e.g., a polyhedron) is equal to the 
so-called  PSD rank of its slack matrix; see the recent 
survey~\cite{fawzi2015positive}. (Note that for non-polyhedral sets, the slack 
matrix is infinite since the number of extreme points and facets is not 
finite.) 

\end{itemize}
These ideas, for example, recently allowed Hamza Fawzi to prove that the PSD 
code 
cannot be represented using the second-order cone~\cite{fawzi2016representing}; 
the proof relies on the fact that the second-order cone rank of the cone of 
3-by-3 PSD matrices is infinite.

\section{Conclusion} 

In this paper, we have introduced the NMF problem and discussed several of its 
aspects. 
The opportunity for meaningful interpretations is 
the main reason why NMF became so popular and has been used in many 
applications. 
NMF is tightly connected with difficult geometric problems; hence 
developing fast and reliable algorithms is a challenge. 
Although important challenges remain to be tackled (e.g., 
developing exact algorithms for small-scale problems), 
even more challenges exist in generalizations of NMF. 
In particular, we mentioned cone factorizations (such as the PSD factorization  
and its symmetric variant~\cite{Gribling2017122}), which are more recent 
problems and have not been explored to their full extent. 

\small

\paragraph{Acknowledgments.}
The author would like to thank the editors of the SIAM Activity Group on Optimization's Views and News (see \url{http://wiki.siam.org/siag-op/index.php/View_and_News}), Stefan Wild and Jennifer Erway, 
for providing insightful feedback that helped improve the presentation of the paper. 
This paper appeared in SIAG/OPT Views and News 25 (1), pp.\@ 7-16. 
The author also acknowledges the support by the F.R.S.-FNRS (incentive grant for scientific 
research n$^\text{o}$ F.4501.16) and by the ERC (starting grant n$^\text{o}$ 
679515).

\bibliographystyle{spmpsci}
\bibliography{bibliographyNMF}

\normalsize

 \end{document}